\documentclass{article}

\usepackage{amsmath,amssymb,enumerate,amsthm}

\usepackage{geometry}
\newtheorem{definition}{Definition}
\newtheorem{remark}{Remark}

\newtheorem{theorem}{Theorem}
\newtheorem{proposition}{Proposition}
\newtheorem{corollary}{Corollary}

\newcommand{\R}{\mathbb{R}}

\newcommand{\T}{\mathbb{T}}
\newcommand{\ep}{\varepsilon}
\newcommand{\Id}{\text{Id}}

\DeclareMathOperator{\diver}{\operatorname{div}}

\DeclareMathOperator{\tr}{\operatorname{Tr}}

\title{A new proof of the uniqueness of the flow  for ordinary differential equations with $BV$ vector
fields}
\author{Maxime Hauray  \and  Claude Le Bris}
\date{March 2010}

\author{
	Maxime Hauray\footnote{Laboratoire d'Analyse, Topologie, Probabilit\'es, CNRS et  Universit{\'e} de 	Provence, Technop\^ole Ch\^ateau-Gombert, 39, rue F. Joliot Curie, 13453 Marseille Cedex 13, 
	France. \tt{hauray@cmi.univ-mrs.fr}  
	}
  	\quad and \quad Claude Le Bris\footnote{CERMICS, {\'E}cole Nationale des Ponts et
	Chauss{\'e}es, 6 \& 8, avenue Blaise Pascal, 77455 Marne-La-Vall{\'e}e Cedex
	2 and INRIA Rocquencourt, MICMAC project, Domaine de Voluceau,  B.P. 105,
 	78153 Le Chesnay Cedex, France. {\tt lebris@cermics.enpc.fr} }
 }

\begin{document}

\maketitle

\begin{abstract} We provide in this article a new proof of  the
uniqueness of the flow solution to  ordinary differential equations with $BV$
vector-fields that have  divergence in $L^\infty$ (or in $L^1$), when
the flow is assumed  nearly incompressible (see the text for the definition of
this term). The novelty of the proof lies in the fact it
does not use the associated transport equation. 
\end{abstract}


\medskip

\section{Introduction and statement of our main result}

 In 1989, P.-L. Lions and R.
DiPerna  showed in \cite{DPL} the existence and the uniqueness of the
almost everywhere defined flow solution to an  ordinary differential equation of the type:
\begin{equation} \label{EDOsimple}
\dot{y}(t) = b(t,y(t)) \, ,
\end{equation}
for $W^{1,1}$ vector fields $b$  with $L^1_{loc}(\R_t,L^\infty_y)$ divergence (along
with some technical
assumptions). For such 'singular' vector fields, the only
possibility is to solve the equation \emph{almost everywhere} on the
space $\Omega$ of initial conditions. In that case, one defines  a flow $X(t,x) : \R \times \Omega \rightarrow \Omega$
satisfying
\begin{equation} \label{EDO}
\left\{ \begin{array}{l}
\dot{X}(t, x) = b(t,X(t,x)) \qquad  \text{for all } t, \\
X(0,x) = x  \qquad
\end{array} \right. .
\end{equation}
for almost all  $x\in\Omega$. An initial time $s\not=0$ may of course
be chosen, and the flow then depends parametrically on this initial time~$s$. With a view to simplifying the presentation, we will assume
henceforth and throughout this article that the field $b$ is
time-independent. Our arguments may be modified to cover the time-dependent case.

In the present article, we also adopt a  notion of
\emph{almost everywhere flow solution} similar to that of DiPerna and
Lions. We denote by $\left(X(t,\cdot)_\# \lambda\right) (E) = \lambda(X(-t,E))$ the pushforward
of the Lebesgue measure $\lambda$. In  the sequel, the
vector-field $b$ will always be assumed at least $L^1_{loc}$.

\begin{definition}[Almost everywhere, nearly incompressible flows]
 \label{def:aeflow}
An almost everywhere flow solution to \eqref{EDO} is a measurable function
$X(t,x) : \R \times \Omega \rightarrow \Omega$ satisfying the following
conditions:
\begin{itemize}
\item[(i)] For almost all $x\in\Omega$, the map $t \mapsto X(t,x)$
is a continuous solution to $\dot{\gamma} = b(\gamma)$ satisfying
  $\gamma(0) =x$ and 
$$ \text{for almost all }x\in\Omega,\quad \; \forall t \in \R, \quad  X(t,x) = x + \int_0^t
b(X(s,x))\,ds $$

\item[(ii)] For all $t$, the measure $X(t,\cdot)_\# \lambda$ is absolutely
continuous with respect to $\lambda$, and, for all $T>0$,  there exists a time-dependent
function $C(T) > 0$ such that the following condition of
\emph{near-incompressibility} holds
$$
\forall T \in \R_+, \; \forall |t| \leq T, \quad \frac{1}{C(T)} \lambda \leq
X(t,\cdot)_\# \lambda \leq C(T) \lambda \, ,
$$
\item[(iii)] $X$ is a one-parameter transformation group, i.e. satisfies
$$
X(t,X(s,x)) = X(s+t,x), \quad \text{for almost all }x\in\Omega, \forall s,\,t
$$

\end{itemize}
\end{definition}

\begin{remark}
Di Perna and Lions originally define in \cite{DPL}  a flow solution with
condition $(i)$ replaced by \emph{$X\in \mathcal{C}(\R,L^1)$
satisfies the ordinary differential equation  \eqref{EDO} in the sense of distribution}. Their definition is
equivalent to ours. It is indeed shown  in \cite{DPL} that the original  definition implies
$(i)$, and it can be shown that conditions $(ii)$ and $(iii)$ together imply that $X \in
\mathcal{C}(\R,L^1)$.
\end{remark}

\begin{remark}
Condition $(ii)$ is in particular satisfied if $b$ is sufficiently smooth, and
$\diver(b) \in L^\infty$. In that case  $C(t) =e^{
\|\diver(b)\|_\infty |t| }$ is convenient.
Besides this case, for a class of ordinary differential equations
coming from  some particular types of hyperbolic equations, such as the
Keyfitz-Krantzer system,  $\diver(b)$ is only $L^1$, but an
estimate of the form $(ii)$ may be established using a
maximum principle. See the work \cite{Bre03} by Bressan for more
details on these systems and that by L.~Ambrosio, F.~Bouchut and C.~De~Lellis
\cite{AmbBouDeL04} for a discussion on the relevance of condition~$(ii)$.
\end{remark}

We now  give a brief state-of-the-art survey on the theory of
ordinary differential equations  with vector fields of low regularity.
The seminal work \cite{DPL} by DiPerna and Lions  has been
complemented and extended notably  by L. Ambrosio
in~\cite{Amb}. Several other authors have made important contributions. We would like
to specifically cite the work~\cite{Ler} by N. Lerner which  has inspired our own, present work.
To date, the minimal conditions that are known to guarantee  the existence and the uniqueness
of the flow are the $BV$ regularity of the vector field, a $L^1$ bound on
the divergence together with a near-incompressibility condition (or more
classically a condition of bounded divergence) of the type $(ii)$.
The classical proofs of such results are based upon the consideration of the
associated transport equation, written either in the conservative form
\begin{equation} \label{Ctransport}
 \frac{\partial u}{\partial t} + \hbox{\rm div}\,( b(x) u) = 0,
\end{equation}
or in the non-conservative form
\begin{equation} \label{NCtransport}
\frac{\partial u}{\partial t} + b(x) \cdot \nabla_x u = 0,
\end{equation}
both with the  initial condition $u(0,x) = u_0(x)$. Remark that, for
divergence-free fields, the two equations coincide. When the
existence and the uniqueness of the solution to the transport equation is
established, for any given initial condition, one deduces the same result for
the a.e. flow solution to the ordinary differential equation. The key ingredient for the resolution
of the transport equation is a commutation lemma (first stated in \cite{DPL}),
which states that
$$
\rho_\ep  * \diver( b u) - \diver( (b * \rho_\ep) u) \underset{\ep\to
0}{\longrightarrow} 0
\;\; \text{in} \; L^1.
$$
C. De Lellis and G. Crippa have recently given in \cite{CriDeL08} a new proof of the
existence and uniqueness of the flow solution of \eqref{EDO}, not using the
 the associated transport equation. Their very interesting
approach provides regularity estimates for $W^{1,p}$ vector-fields with
$p>1$ but  seemingly fails  for $W^{1,1}$ vector-fields,  unfortunately.
Very recently, the approach has been succesfully improved by
P.-E.~Jabin. In \cite{Jab09}, this author extends the direct method by C. De Lellis and G. Crippa to the case of bounded SBV vector fields (with locally finite jump set for the $d-1$ Hausdorff measure) in any
dimension, and also to two-dimensional $BV$ vector-fields that satisfy a
particular, local assumption in the direction of the flow. In addition,
in \cite{ChaJab09},  the same
author, in collaboration with N.~Champagnat, has proved that, in the particular case of an
ODE corresponding to the Newton equation of motion (that is, $\dot X=V$,
$\dot V = F(X)$), there exists a flow
solution if $F \in H^{3/4}\cap L^\infty$, and this flow is unique among
the class of flows obtained by regularization of the problem.

\medskip

Before we get to the heart of the matter, some comments are in order
regarding our assumptions, which are slightly different from those of~\cite{CriDeL08}. The important assumption in \cite{CriDeL08} is the
bound from above for $t>0$ in (ii), in order to prevent \emph{contraction} of
the measure.
In addition, the group structure (our assumption (iii)) is
not required in~\cite{CriDeL08}. In the present state of our
understanding, such differences seem to us related to the different techniques
of proof used in~\cite{CriDeL08} and here.
The need for (iii) is consistent with our previous work~\cite{HauLeBLio}
where, in order for the proof of uniqueness to be performed,  the group
structure is necessary in the
Lagrangian viewpoint, but not in the Eulerian viewpoint. The assumption on the group structure is also present in
\cite{DPL}. In the present contribution, we need this assumption because
we transform integrals  back and forth in time
(see our argument in Section 2 below).
To end this discussion, we also note that Assumption~(ii) may be
slightly weakened  if we are only interested by
the uniqueness of the forward-in-time flow. In that case, we may replace (ii)
by
$$ (ii') \qquad \frac{1}{C(T)} \lambda \leq X(t,\cdot)_\# \lambda <<
\lambda , \qquad \forall t \in [0,T] \,.$$
This is in contrast to the assumptions of~\cite{CriDeL08}, which
uses a  bound from above: $X(t,\cdot)_\# \lambda \leq C(T)\,\lambda$.  This owes
to the very
strategy of our proof, which mostly uses the backward-in-time flow.
For instance, the most important step \eqref{eq:impstep} is performed on
the backward flow. In any event, owing to the group property, the upper and lower bounds
are closely related to one another. The two uniform
bounds in $(ii)$ are indeed necessary to conclude when the issue considered
is the uniqueness of both the
backward-in-time and the forward-in-time flows.  We refer the reader to Remark~\ref{rk:final} at the end of the proof for more details.

\subsection{Main result}
  The purpose of this article is to  give a new and direct proof of the
uniqueness of the a.e. flow solution to~\eqref{EDO} for $BV$  vector fields,
\emph{without} arguing on  the associated transport equation. We adopt
the approach  already used in \cite{HauLeBLio} for $W^{1,1}$ vector fields. Basically, the
 commutation lemma instrumental in the proof contained in  this prior
 publication is replaced by another strategy of proof, namely the introduction of a
second variable. This is explained in details in the next paragraph.

\medskip

Our result is the following:

\begin{theorem}
\label{thm:1}
Let $b$ be a $BV$ vector field on the $N$-dimensional torus $\T^N$. If $\diver(b) \in L^1$, then there
exists at most one a.e. flow solution to~\eqref{EDO}, in the sense of
Definition~\ref{def:aeflow}.
\end{theorem}

\begin{remark}
In
\cite{CapPer}, B. Perthame and I. Capuzzo Dolcetta remarked that the
assumption  "$b \in W^{1,1}$" of the original work by DiPerna and Lions
could be replaced by the weaker assumption "the
symmetric part of $Db$ is a matrix-valued $L^1$ function". This
observation seems to not be valid for the present strategy of proof, and more
generally in the $BV$ case. The reason is, their argument is based on the use of
radially symmetric regularization kernels, while the regularization
kernels we use here for the BV case are typically anisotropic. 
\end{remark}

\begin{remark}
  Theorem~\ref{thm:1} states uniqueness of the flow. We do not know
  which condition, besides the general assumptions of
  Theorem~\ref{thm:1}, makes possible the proof of existence of a flow
  in the sense of Definition~\ref{def:aeflow}, using the approach
  developed here for uniqueness.
\end{remark}

\subsection{Main idea of the proof}
  To start with, we outline here the proof performed in details in the
  next section. As already said, the proof uses a technique introduced
in \cite{HauLeBLio}. In that work, a smooth
convolution kernel $\rho$, with normalized integral, is considered. It
is then proved
that for any two a.e. flows $X$ and $Y$ solutions to the ordinary differential equation with $W^{1,1}$
coefficients,
$$
\lim_{\ep \rightarrow 0} \frac{d}{dt} \left( \int \int | X(t,x) - Y(t,y)|
\frac{1}{\ep^N} \rho(\frac{x-y}{\ep}) \,dx \,dy \right) = 0.
$$
Now, the limit of the integral is
$$
\frac{d}{dt} \left( \int | X(t,x) - Y(t,x)| \,dx \right) = 0.
$$
This shows that, for all $t$,
$$
\int | X(t,x) - Y(t,x)|\,dx =0,
$$
since this quantity vanishes at initial time. The uniqueness of the solution
follows. Remark that the introduction of the extra-variable $y$ allows to
perform the calculation  without using the
transport equation.

\medskip

Our aim is to now modify  the above approach and treat BV vector fields. For this
purpose, we use a convolution kernel well adapted to the geometry of the flow
and the possible singularities of the BV vector field under
consideration.
In short, we consider the
regularization kernel
$$
 \frac{1}{\ep^N} \rho(x,\frac{x-y}{\ep}) \text{ with }\
\rho(x,z) = F_0(|U(x) z|^2) \det U(x) \, , \text{ and }
U(x) = Id + \gamma \eta(x) \otimes \eta(x) \, .
$$
Here, $F_0$ is a smooth function, $\gamma$ is a constant that will  be sent to
infinity, and $\eta$ is an approximation of the direction normal to the jumps of
the measure $Db$. The purpose of such a construction is to have a regularization
that decreases faster in the direction normal to the jumps. The idea of a
direction-dependent regularization was first introduced by P.L.~Lions in
\cite{Lio}. N. Lerner introduced the specific position-dependent regularization
used here in \cite{Ler} with a view to simplifying the proof of uniqueness
originally given by L.~Ambrosio for the $BV$ case. His argument,
however, is  still based upon the
equivalence with the transport equation. In the present paper, we combine his
argument with the approach consisting in introducing a second variable,
already employed in \cite{HauLeBLio} for $W^{1,1}$
vector fields.

\section{Proof of Theorem~\ref{thm:1}}

This section is devoted to the proof of Theorem \ref{thm:1}. We denote by
$\mu_1(t,\cdot)$ (resp. $\mu_2(t,\cdot)$) the $L^\infty$ density of
the measure $X(-t,\cdot)_\# \lambda$ (resp. $Y(-t,\cdot)_\# \lambda$) with respect to
$\lambda$.
\medskip

Consider now the kernel
$$
\frac{1}{\ep^N} \rho(x,\frac{x-y}{\ep}),
$$
where $\rho$ is a smooth, compactly supported,  function, from $\T^N \times
\T^N$ to $\R^+$ which we will make precise below. Assume in addition $\rho$
satisfies $\int \rho(x,z) \,dz =1$ for all $x$. Our aim is to estimate
\begin{equation} \label{eq:defIep}
I_\ep (t) = \frac{d}{dt} \left( \int \int | X(t,x) - Y(t,y)| \, \frac{1}{\ep^N}
\rho(x,\frac{x-y}{\ep}) \mu_1(t,x)  \mu_2(t,y) \,dx \,dy \right).
\end{equation}
where $X$ and $Y$ are two flow solutions to~\eqref{EDO}. In the sense of distributions,
\begin{equation} \label{eq:limit}
\lim_{\ep \rightarrow 0} I_\ep (t) = \frac{d}{dt} \left( \int | X(t,x) -
Y(t,x)| \mu_1(t,x) \mu_2(t,x) \, dx \right).
\end{equation}
This is established using the Lebesgue continuity of the
functions $Y$ and $\mu_2$ at almost every point,  along with  the $L^\infty$ bound on $\mu_1$. Remark that the
Lebesgue continuity may be used if the support of $\rho(x,\cdot)$ is not
exceedingly stretched in one direction (more specifically, we should
have some constant $c >
0$ such that  $\forall x\in\T^N, \;
B(0,c^{-1}) \subset \mathrm{Supp} \rho(x,\cdot) \subset B(0,c)$, See \cite{Stein70} for more details). The kernel we shall use satisfies such a  condition for all $\varepsilon>0$,
even though in the limit of a vanishing $\varepsilon$, it is infinitely stretched.
Our purpose is to show that the limit \eqref{eq:limit} is
$$
\lim_{\ep \rightarrow 0} I_\ep (t) = - \int |X(t,x)
-Y(t,x)| \diver(b)(x) \mu_1(t,x)\mu_2(t,x) \,dx.
$$
This will eventually prove the uniqueness of the flow solution to \eqref{EDO} using the bounds from
below on  $\mu_1$ and $\mu_2$ inferred from  $(ii)$.
To this end, we first perform the change of variable  $(x,y) \rightarrow (X(t,x),Y(t,y))$
in $I_\ep(t)$, and then differentiate under the integral
$$
I_\ep(t)  =  \frac{d}{dt} \left( \int \int | x - y| \, \frac{1}{\ep^N}
\rho(X(-t,x),\frac{X(-t,x)-Y(-t,y)}{\ep}) \,dx \,dy \right),
$$
which we write $I_\ep(t)  =  I_\ep^1(t) + I_\ep^2 (t)$, with
\begin{eqnarray*}
 \label{eq:impstep}
I_\ep^1(t) & = &   - \int \int  \frac{| x - y|}{\ep^N} \partial_1
\rho(X_{-t}(x), \frac{X_{-t}(x)-Y_{-t}(y)}{\ep}) \cdot b(X_{-t}(x)) \,dx \,dy
\\
I_\ep^2(t) & = &  - \int \int  \frac{| x - y|}{\ep^{N+1}} \partial_2 \rho
(X_{-t}(x), \frac{X_{-t}(x)-Y_{-t}(y)}{\ep}) \cdot (b(X_{-t}(x)) - b(Y_{-t}(y))) \,dx \,dy \,
\end{eqnarray*}
with the notations $X_t(x) =X(t,x)$ and $Y_t(y)=Y(t,y)$.
Then, we return to the original variables $(x,y)$
\begin{eqnarray*}
I_\ep^1(t) & = &  - \int \int | X_t(x) - Y_t(y)| \, \frac{1}{\ep^N} \partial_1
\rho(x, \frac{x-y}{\ep}) \cdot b(x) \mu_1(t,x)  \mu_2(t,y) \,dx \,dy \\
I_\ep^2(t) & = & - \int \int | X_t(x) - Y(t,y)| \, \frac{1}{\ep^{N+1}}
\partial_2 \rho (x , \frac{x-y}{\ep}) \cdot (b(x) - b(y))  \mu_1(t,x)
\mu_2(t,y) \,dx \,dy ,
\end{eqnarray*}
and next use the change of variable $z = (y-x)/\ep$
\begin{eqnarray}
I_\ep^1(t) & = & - \int \int | X_t(x) - Y_t(x + \ep z)| \, \partial_1 \rho(x,
z) \cdot b(x)  \mu_1(t,x)  \mu_2(t,x + \ep z ) \,dx \,dz \label{eq:I1}\\
I_\ep^2(t) & = & - \int \int | X_t(x) - Y_t(x+ \ep z)| \, \partial_2 \rho (x ,
z) \cdot \frac{(b(x) - b(x + \ep z))}{\ep}  \dots \nonumber \\ 
&& \hskip 55truemm \dots \mu_1(t,x)  \mu_2(t,x + \ep z )
\,dx \,dz \,. \label{eq:I2}
\end{eqnarray}
We now need to estimate these two terms when $\ep$ goes to zero. We begin with
the easiest of the two, namely~$I_\ep^1$.

\medskip

\noindent \textbf{Step 1: Limit of $I_\ep^1$}

Because $\rho$ is smooth,  $b \in L^1$ , and almost all points are Lesbesgue
points for the two  functions $Y$ and $\mu_2$, we can use the Lebesgue dominated
convergence theorem and obtain
$$
\lim_{\ep \rightarrow 0} I_\ep^1 (t) =  - \int | X_t(x) - Y_t(x)| \left( \int
\partial_1 \rho(x, z) \,dz \right) \cdot b(x) \mu_1(t,x)\mu_2(t,x)\,dx.
$$
Now
$$
\displaystyle \int \partial_1 \rho(x, z) \,dz = \frac{d}{dx} \left( \int
\rho(x, z) \,dz \right) = 0 \, ,
$$
since $\int \rho(x, z) \,dz = 1$, for all $x$. Thus,
\begin{equation}
  \label{eq:Ieps1}
  \lim_{\ep \rightarrow 0} I_\ep^1 (t) = 0.
\end{equation}

The treatment for $I_\ep^2$ is more elaborate and will necessitate several
steps.

\medskip
\noindent\textbf{Step 2: Bound for $I_\ep^2$}

We now wish to pass to the limit $\ep \rightarrow 0$ in \eqref{eq:I2}. If $b$
were $W^{1,1}$, the limit could easily be identified. It would suffice to replace
$(b(x + \ep z) - b(x))/\ep$ by  $\int_0^1 Db(x + \theta \ep z) \cdot z \,d\theta
$ in \eqref{eq:I2}, and next use the Lebesgue dominated convergence theorem. All
this does not require making  specific the convolution kernel $\rho$ (See below
and \cite{HauLeBLio}). Owing to the presence of the singular part of $Db$, we
have to argue more carefully.

\medskip

To proceed further, we  recall the following result:
\begin{proposition}{\bf [from \cite[Theorem~1.28, Corollary~1.29]{AFP}]}
\label{prop:radon}
Let $b$ be a $BV$ vector-field on $\T^N$.

\noindent (i) The Radon-Nikodym decomposition of its derivative $Db$ writes
$$
Db = D^a b + D^s b, \quad \text{with} \quad D^a b << \mathcal{L}^N, \; D^s b
\perp  \mathcal{L}^N ,
$$
where the superscript $a$ stand for "absolute continuous part", and $s$ stand
for "singular" respectively. As $D^a b$ is absolutely continuous with respect to
the Lebesgue measure, we write it
$$
D^a b = \partial^a b \, dx \, ,
$$
where $\partial^a b$ is a $L^1$ matrix-valued fonction.

\noindent (ii)  In addition, the polar decomposition of the singular part $D^s
b$ of the measure $Db$ writes
$$
D^s b = M^s \,|D^s b| \, ,
$$
where $|D^s b|$ is the total variation of the matrix-valued measure $D^s b$,
and $M^s$ a matrix-valued fonction, such that $|M^s(x)| = 1$, $|D^s b|$-a.e (the
norm used for $M^s$ is the norm induced on matrices by the Euclidian norm of
$\R^n$).
\end{proposition}

In view of the above decomposition, we now claim that
\begin{equation}  \label{eq:variante}  \begin{split}
\limsup_{\varepsilon\longrightarrow 0}    &  \int \! \int  \int_0^1  |
X_t(x)  - Y_t(x+ \ep z)|  \dots  \\
\dots & \left|  \partial_2 \rho (x , z) \cdot \left(
\frac{b(x + \ep z) - b(x)}{\ep}  -    \partial^a b(x + \ep \theta z) \! \cdot \! z
\right) \right| 
 \mu_1(t,x)  \mu_2(t,x + \ep z ) \,d\theta \,dx \,dz \\
& \hskip 4truecm \leq 2 C(t)^2 \,\int \! \int \left| \partial_2 \rho (x , z) \cdot M^s(x)
\cdot z \right| \,d|D^s b|(x) \,dz .
\end{split}  \end{equation}
For convenience, we denote by
$$
I^2_{\varepsilon,a}=-\int \! \int \int_0^1  |
X_t(x) - Y_t(x+ \ep z)| \,   \partial_2 \rho (x , z) \cdot
\partial^a b(x + \ep \theta z) \cdot z\,d\theta \,dx \,dz,
$$
in the left-hand side, and
$$\bar{I}_s^2(t)=\int \! \int \left| \partial_2 \rho (x , z) \cdot M^s(x)
\cdot z \right| \,d|D^s b|(x) \,dz ,
$$
in the right-hand side.

\medskip

To prove our claim, we regularize $X$, $Y$, $\mu_1$ and $\mu_2$, using some
smooth $X^\alpha$,
$Y^\alpha$, $\mu^\alpha_1$ and $\mu^\alpha_2$. Next, we replace $(b(x+
\ep z) - b(x))/\ep$ by  $\int_0^1 Db(x + \theta \ep z) \cdot z \,d\theta $
 (an equality
true for almost all $(x,z)$) and perform the change of variable $x' = x + \ep
\theta z$ (we use it even for the measure $Db$ because this is a linear change
of variable). We obtain 
\begin{eqnarray} \label{I2ep}
I_\ep^{2,\alpha}(t) & := & - \int \! \int | X^\alpha_t(x) - Y^\alpha_t(x+
\ep z)|\, \partial_2 \rho
(x , z) \cdot \frac{(b(x + \ep z) - b(x))}{\ep}  \dots \nonumber \\ 
&& \hskip 60truemm \dots \mu_1^\alpha(t,x)
\mu_2^\alpha(t,x + \ep z ) \,dx \,dz \nonumber\\
 & = &  - \int \! \int \!\int_0^1 | X^\alpha_t(x) - Y^\alpha_t(x + \ep z)| \,
\partial_2 \rho (x , z) \cdot  Db(x + \theta \ep z) \cdot z  \dots \nonumber \\ 
&& \hskip 60truemm \dots \mu^\alpha_1(t,x)
\mu^\alpha_2(t,x + \ep z )\,dx \,dz \,d\theta \nonumber\\
& = & - \int \!\int \!\int_0^1 | X^\alpha_t(x - \ep \theta z) -
Y^\alpha_t(x + \ep (1-\theta) z)| \dots \nonumber\\
&&\hskip -3truemm \dots \partial_2 \rho (x - \ep \theta z, z) \cdot
Db(x) \cdot z  \mu^\alpha_1(t,x - \ep \theta z) \mu^\alpha_2(t,x + \ep
(1-\theta) z )\,dx \,dz \,d\theta \, .
\end{eqnarray}
Let us decompose $I_\ep^{2,\alpha}$ in two parts, according to the above
Proposition~\ref{prop:radon},
$$
I_\ep^{2,\alpha}(t) = I_{\ep,a}^{2,\alpha}(t) + I_{\ep,s}^{2,\alpha}(t),
$$
where
\begin{eqnarray}
I_{\ep,a}^{2,\alpha}(t) & = &  - \int \!\int \!\int_0^1 |
X^\alpha(t,x - \ep \theta z) - Y^\alpha(t,x + \ep (1-\theta) z)| \dots \nonumber\\
&& \hskip -7mm \dots \partial_2
\rho (x \!- \!\ep \theta z, z) \!\cdot\!  \partial^a b(x) \! \cdot \! z \, \mu^\alpha_1(t,x \!- \!
\ep \theta z) \mu^\alpha_2(t,x \!+\! \ep (1\!-\!\theta) z )\,dx \,dz \,d\theta
\nonumber\\
|I_{\ep,s}^{2,\alpha}(t)| & \leq &  2 C(t)^2 \int \!\int \!\int_0^1 \partial_2
\rho (x - \ep \theta z, z) \cdot  M^s(x) \cdot z \, |D^s b|(x)\,dz \,d\theta \,
\nonumber \end{eqnarray}
where we have used that $|X^\alpha-Y^\alpha| \leq 2$ (as we work on
the torus).  And letting $\ep$ going to zero, we obtain (\ref{eq:variante}) for $X^\alpha$, $Y^\alpha$ and
the $\mu_i^\alpha$. Then (\ref{eq:variante}) is obtained letting $X^\alpha$,
$Y^\alpha$, $\mu_1^\alpha$ and $\mu_2^\alpha$ approximate $X$, $Y$,  $\mu_1$ and
$\mu_2$, respectively.

\medskip

The majoration~(\ref{eq:variante}) being established, we proceed as follows.
Arguing as above for $I^1_\varepsilon$, that is using the smoothness of
$\rho$ and the fact that almost every point is a Lebesgue point for $Y$,
$\mu_1$ and
$\mu_2$, we obtain
\begin{eqnarray*}
\lim_{\ep \rightarrow 0}I^2_{\varepsilon,a}(t) & = & \lim_{\ep \rightarrow 0}
- \int \! \int \! \int_0^1 |X_t(x - \ep \theta z) - Y_t(x + \ep
(1-\theta) z)| \dots \nonumber\\
&& \dots \partial_2 \rho (x -\varepsilon\theta z, z) \cdot \partial^a
b(x)\cdot z \, \mu_1(t,x - \ep \theta )  \mu_2(t,x +  \ep (1-\theta) z ) \,d\theta
\,dx \,dz\\
&= & - \int \! \int |X_t(x) - Y_t(x)| \partial_2 \rho(x,z) \cdot \partial^a b(x)
\cdot
z \, \mu_1(t,x)\mu_2(t,x) \,dx \,dz \\
& = & - \int |X_t(x) - Y_t(x)|  R_a(x)  \mu_1(t,x) \mu_2(t,x) \,dx,
\end{eqnarray*}
with $ R_a(x) = \int \partial_2 \rho(x,z) \cdot \partial^a b(x) \cdot z \,dz$.
To calculate this term, we integrate by parts and use the property $\forall x
\in \T^N, \; \int \rho(x,z)  \,dz =1$ 
\begin{eqnarray}
R_a(x) & = & \sum_{i,j} \int \frac{\partial \rho}{\partial z_i}(x,z)
\frac{\partial^a b_i}{\partial x_j}(x) z_j \,dz =  \sum_{i,j} \frac{\partial^a
b_i}{\partial x_i}(x) \int - \rho (x,z) \frac{\partial z_j}{\partial z_i} \,dz
\nonumber  \\
& = & 
 - \sum_i \frac{\partial^a b_i}{\partial z_i}(x) \int \rho(x,z)  \,dz  = - \diver^a b. \label{eq:Ra}
\end{eqnarray}
So we have obtained
\begin{equation}
  \label{eq:Iepa}
  \lim_{\ep \rightarrow 0} I_{\ep,a}^2 =  \int |X_t(x) - Y_t(x)|  \diver^a b(x)
\mu_1(t,x) \mu_2(t,x) \,dx.
\end{equation}
The next step consists in proving that the right-hand side of
(\ref{eq:variante}) may be chosen arbitrarily small.

\medskip
\noindent\textbf{Step 3: A bound on the singular part}

In order to estimate the right hand side of (\ref{eq:variante}), we now use a geometric
information, namely, the special form of $M^s(x)$, proved by G. Alberti
\cite{Alb}.
\begin{theorem} {\bf [Alberti's rank one  Theorem, \cite[Theorem 3.94]{AFP}]}
Let $b$ be a $BV$ vector-field defined on $\T^N$, and write $Db = D^s b + D^a
b$ the Radon-Nikodym decomposition of its gradient. Consider $D^s b = M^s \, |D^s
b|$ the polar decomposition of the singular part as in
Proposition~\ref{prop:radon}. Then, $M^s$ is of rank one $|D^s b|$-almost
everywhere, that is,  there exists two vector-valued functions $\xi_b$ and
$\eta_b$, both $|D^s b|$-measurables, such that $\xi_b$ and $\eta_b$ are unit vectors
$|D^s b|$-a.e. and satisfy
$$
M^s(x) = \xi_b(x) \otimes \eta_b(x), \quad |D^s b|-\text{almost everywhere},
$$
where $\xi_b \otimes \eta_b$ denotes the linear map $z \mapsto \langle \eta_b,
z \rangle \xi_b $.
\end{theorem}
\begin{corollary} \label{divsing}
As a consequence, the singular part of the divergence is 
$$\diver^s b = \langle \xi, \eta \rangle |D^s b| \,.$$
If we assume that the divergence of $b$  belongs to
$L^1$, it follows that
$$
\langle \xi ,\eta \rangle =0, \quad |D^s b|-\text{almost everywhere},
$$
a property that  will be crucial in the sequel.
\end{corollary}

Using the decomposition provided by Theorem 2, we rewrite our bound in
\eqref{eq:variante}, which we denote $\bar{I}_s^2(t)$ in the sequel
$$
\bar{I}_s^2 (t) \leq  2 C(t)^2 \int \int | \langle \partial_2 \rho (x,z),
\xi_b(x) \rangle | | \langle \eta_b(x), z \rangle | \,d|D^s b|(x) \,dz.
$$
In order to render the right-hand side arbitrarily small, we now make specific our
convolution kernel $\rho$. We choose
$$
\rho(x,z) = F_0(|U(x)z|^2) \det(U(x)),
$$
where $F_0$ is a smooth, compactly supported, non negative  function
such that \\ \mbox{$\int_{\R^N} F_0(|z|^2) \, dz = 1$}, and $U$ is a smooth,
matrix-valued function, such that $U(x)$ is an  orientation preserving
matrix for all $x$. Note that owing to the presence of the determinant,
the integral of $\rho(x,\cdot)$ remains equals to one independently of
$x$. The dilation matrix $U(x)$ is set to $U(x) = Id + \gamma \eta(x)
\otimes \eta(x) $ (with the notation $a \otimes b$ for the endomorphism
$x \rightarrow \langle b,x \rangle a$), where $\eta$ is a smooth
vector-valued function. On the jump of the measure $Db$, $\eta$ will be
chosen later as an approximation of the direction normal to the jump
set. The factor $\gamma$ will be chosen as large as possible. It may possibly
depend upon $x$ and be large only on a neighbourhood of the singular set of the
measure $Db$, but we for simplicity of the calculation we will not use that not essential
possibility here.

The partial derivative of $\rho$ writes
$$
\partial_2 \rho (x,z) = 2 F_0'(|U(x)z|^2) \langle U(x)z, U(x) \cdot \rangle
\det(U(x)).
$$

We use this in the bound on $\bar{I}_s^2(t)$ to obtain
\begin{equation*}
\bar{I}_s^2 (t) \leq  C \int \int  |F_0'(|U(x)z|^2)| | \langle U(x)z , U(x)
\xi_b(x) \rangle | \langle \eta_b(x), z \rangle |\det(U(x)) \,d|D^s b|(x) \,dz,
\end{equation*}
where here and below $C$ denotes various irrelevant constants.
To simplify this term, we perform the change of variable $z \rightarrow U(x)
z$, and obtain
\begin{equation} \label{boundetaxi}
\bar{I}_s^2 (t) \leq C \int \int  |F_0'(|z|^2)| | \langle z , U(x) \xi_b(x)
\rangle | \, | \langle \eta_b(x), U^{-1} z \rangle | \,d|D^s b|(x) \,dz.
\end{equation}

We next intend to use  the special form $U(x) = Id + \gamma(x) \eta(x)
\otimes \eta(x)$, to bound from above the two scalar products.
Let us first formally illustrate our argument, performing our
calculation with $\eta=\eta_b$, as if $\eta_b$ were smooth. In this case,
$$ | \langle z , U(x) \xi_b(x) \rangle | = |\langle z , \xi_b(x) \rangle | \leq
|z| ,$$
because $\langle \eta_b, \xi_b\rangle = 0$ and $\xi_b$ has unit norm. For the
second scalar product,
$$ | \langle \eta_b(x), U^{-1}(x) z \rangle | = \frac{1}{1 + \gamma} | \langle
\eta_b(x),  z \rangle | \leq \frac{1}{1 + \gamma}, $$
because $ \displaystyle U^{-1} = \Id - \frac{\gamma}{1+\gamma}\eta_b
\otimes \eta_b $.

Inserting  these bounds in \eqref{boundetaxi}, we obtain
$$
\bar{I}_s^2 (t) \leq \frac{C}{1 + \gamma} \int \int  |F_0'(|z|^2)|
\,|D^s b|(x) \,dz \leq \frac{C(F_0,b)}{1 + \gamma},
$$
where the constant $C(F_0,b)$ depends only of $F_0$ and $b$. It remains
then to let $\gamma$ to infinity to obtain $\bar{I}_s^2 (t)=0$ and
conclude our (formal) proof.

\medskip

We now modify the above formal argument using an approximation  $\eta$
of $\eta_b$, instead of $\eta_b$ itself. First, we remark
\begin{eqnarray}
| \langle  z, U(x) \xi_b \rangle | & = & | \langle  z, \xi_b + \gamma \langle
\xi_b , \eta \rangle \eta \rangle | \nonumber \\
& \leq & (1 + \gamma |\langle \xi_b , \eta \rangle |) \, |z| \nonumber \\
& \leq & (1 + \gamma |\langle \xi_b , \eta - \eta_b \rangle |) \, |z| \nonumber
\\
& \leq & (1 + \gamma |\eta - \eta_b|) \, | z | \label{scalar1},
\end{eqnarray}
where we have used $\langle \xi_b , \eta_b \rangle =0$, $|D^s b|$-a.e (from
Corollary~\ref{divsing}), and that $\xi_b$, $\eta_b$ are unit vectors. To bound
the scalar product $| \langle \eta_b(x), U^{-1} z \rangle |$, we decompose $z$
in $z = z_{\eta} + z_\perp$, where $z_{\eta}$ is the projection of $z$ on $\R
\, \eta$
\begin{eqnarray}
| \langle \eta_b, U^{-1} z \rangle | & = & \left| \left\langle \eta_b, z  -
\frac{\gamma}{1+\gamma} \langle \eta , z \rangle \eta  \right\rangle \right|
\nonumber \\
& = & \left|  \left\langle \eta_b, z_\perp  + \frac{1}{1+\gamma} z_{\eta}
\right \rangle \right|  \nonumber \\
& \leq & \left| \left\langle \eta_b - \eta, z_\perp  + \frac{1}{1+\gamma}
z_{\eta}  \right\rangle \right| + \left| \left\langle \eta, z_\perp  +
\frac{1}{1+\gamma} z_{\eta}  \right\rangle \right| \nonumber \\
& \leq & \left( | \eta_b - \eta | + \frac{1}{1+\gamma} \right) | z |.
\label{scalar2}
\end{eqnarray}

From \eqref{scalar1} and \eqref{scalar2} we deduce
$$
| \langle  z, U(x) \xi_b \rangle | \, | \langle \eta_b, U^{-1} z \rangle |
\leq  \left( 2 | \eta - \eta_b| + \frac{1}{1+\gamma} + \gamma |\eta - \eta_b|^2
\right) |z|^2.
$$

We insert this bound in \eqref{boundetaxi} and obtain
\begin{eqnarray}
\bar{I}_s^2 (t) &\leq & C \int \left( | \eta - \eta_b| + \frac{1}{1+\gamma} +
\gamma |\eta - \eta_b|^2 \right)  \left (\int  F_0'(|z|^2) |z|^2 \,dz.\right)
\,|D^s b|(x) \nonumber \\
&\leq&  C(F_0) \left( \frac{1}{1+\gamma} + (1 +2 \gamma) \int | \eta - \eta_b|
\,|D^s b|(x) \right), \label{lastbound}
\end{eqnarray}
because the integral $\int  F_0'(|z|^2) |z|^2 \,dz$ is fixed,  and  both
$\eta_b$ and $\eta$ are unit vectors.

We finally show that
\begin{equation}
  \label{lastbound0}
  \inf_{\gamma>0,\eta \, \text{ smooth}} \left( \frac{1}{1+\gamma} + (1 +2
\gamma) \int   | \eta - \eta_b| \,|D^s b|(x) \right) =0  .
\end{equation}
To this end, we first choose $\gamma$ such that $1/(1+\gamma)$ is small,
and then construct a smooth function $\eta$, sufficiently close to $\eta_b$ on
the support of $D^s b$ so that $(1 +2 \gamma) \int   | \eta - \eta_b|
\,|D^s b|(x)$ is also small (use for that classical approximation theorem with
respect to the Radon measure $D^s b$).
Note that $\eta$ can be arbitrarily extended to the whole torus as its value
outside the support of $D^s b$ is irrelevant. This concludes the proof
of the convergence of the right-hand side of (\ref{eq:variante}) to zero.

\medskip
\noindent\textbf{Step 4: Conclusion}
Collecting all the previous results, we  obtain
\begin{equation} \label{eq:eqfin}
\begin{split}
\frac{d}{dt} \int |X_t(x) -Y_t(x)| & \mu_1(t,x)\mu_2(t,x) \,dx = \\ 
& - \int |X_t(x) -Y_t(x)| \diver(b)(x) \mu_1(t,x)\mu_2(t,x) \,dx.
\end{split}
\end{equation}
where we have replaced $\diver^a(b)$ by $\diver(b)$, since we are dealing
with vector fields $b$ having at least, divergence in $L^1$.
If $\diver(b) \in L^\infty$, then
$$
\frac{d}{dt} \int |X_t(x) -Y_t(x)| \mu_1(t,x)\mu_2(t,x) \,dx\leq C \int |X_t(x)
-Y_t(x)|  \mu_1(t,x)\mu_2(t,x) \,dx \, .
$$
Since the integral in the right hand side vanishes initially, we conclude that
$$
\int |X_t(x) -Y_t(x)|  \mu_1(t,x)\mu_2(t,x) \,dx =0
$$
and finally that $X(t,\cdot) = Y(t,\cdot)$ a.e. in $x$ since the $\mu_i$ are bounded away from $0$.
 Note that, as usual, if only  the solution at positive times if of interest,
an assumption  on the  negative part $\diver(b)^-$ of the divergence
suffices to conclude.

\medskip

When only the weaker  hypothesis~$\diver(b)\in L^1$ holds, we have to slightly adapt the
above argument. We choose a smooth compactly supported function
$\phi(x)$, insert a factor $\phi(X(t,x))$ in the integral \eqref{eq:defIep} defining
$I_\ep$. We  now estimate
$$
I_\ep^\phi (t) = \frac{d}{dt} \left( \int \int \phi(X_t(x))| X_t(x)  - Y_t(y)| \,
\frac{1}{\ep^N} \rho(x,\frac{x-y}{\ep}) \mu_1(t,x)  \mu_2(t,y) \,dx \,dy
\right).
$$
The above argument carries over to the present case. An  equality
similar to \eqref{eq:eqfin} is obtained
\begin{equation} \label{eq:eqfinL1}
\begin{split}
\frac{d}{dt} \int & \phi(X_t(x)) |X_t(x)  -Y_t(x)| \mu_1(t,x)\mu_2(t,x) \,dx  \\
 & = - \int \phi(X_t(x)) |X_t(x) -Y_t(x)| \diver(b)(x) \mu_1(t,x)\mu_2(t,x) \,dx
\end{split}
\end{equation}
which can also be written (using the change of variable $x=X_t(x)$)
\begin{equation}
\label{eq:eqfinnL1}
\begin{split}
\frac{d}{dt}  \int  & \phi(x) |x - Y_t(X_{-t}(x))| \mu_2(t,X_{-t}(x)) \,dx  \\
&   =  -\int \phi(x) |x -Y_t(X_{-t}(x))| \diver(b)(X_{-t}(x)) \mu_2(t,X_{-t}(x)) \,dx.
\end{split}
\end{equation}
We next define $u(t,x)= |x - Y_t(X_{-t}(x))|
\mu_2(t,X_{-t}(x))$. Equation \eqref{eq:eqfinnL1} holding for all $\phi$,
it follows that
\begin{equation}
\label{eq:no-deriv}
\frac{\partial u}{\partial t} + \diver(b)(X_{-t}(x))\, u= 0,
\end{equation}
in the distributional sense.
There is no derivative of $u$ with respect to $x$ in the equation, so that the variable $x$ is only a parameter.
Since $\diver(b) \in L^1$ and condition (ii) holds, we have
$$\int_x \int_0^T |\diver(b)(X(-t,x))| \,dtdx < +\infty$$
for all time $T$. So that, for almost all $x$,  $\int_0^T |\diver(b)(X(-t,x))| \,dt < +\infty$. Therefore equation (\ref{eq:no-deriv}) is well-posed for
almost all $x$, and since  by construction its solution  $u$ vanishes at
initial time, it vanishes for all time: $u(t,x)=0$ for all $t$, a.e. in
$x$. This concludes the proof: $X\equiv Y$.

\begin{remark}
  \label{rk:final}
Given the above argument, the reader may now understand that, when  the
uniqueness of only the flow for positive times  is under study, the upper
bound $X(t,\cdot)_\# \lambda \leq C(T) \lambda$, for times $t>0$ in (ii) of Definition~\ref{def:aeflow} may be somehow relaxed.
In that case, and as briefly announced in the introduction, we may only assume that
$$ (ii') \qquad \frac{1}{C(T)} \lambda \leq X(t, \cdot)_\# \lambda <<
\lambda , \qquad \forall t \in [0,T] \,. $$ where the symbol $<<$ means
here \emph{absolutely continuous with respect to}.
The group property allows to equivalently state (composing the
previous inequality with $X(- t,\cdot)$) that
$$ (ii') \qquad \lambda << X(-t, \cdot)_\# \lambda \leq C(T) \lambda, \qquad \forall t
\in [0,T] \,. $$
 In that case, we have
$$ 0 < \mu_2(t,x) \leq C(T) \; , \qquad \hbox{\it a.e. in}\quad x \, , $$ and also $0 <
\mu_2(t,X(-t,x))$ almost everywhere  in $x$ since,  considering 
$$(ii') \,,   \qquad \lambda(X(t,A)) \leq C(T) \lambda(A)$$
for all mesurable set $A$. This suffices to show
that, for all $t>0$, and  almost everywhere  in $x$,  $u(t,x)=0$ and thus conclude
$X(t,x)=Y(t,x)$ in the above argument.
\end{remark}

\begin{remark}
As pointed out by L. Ambrosio in \cite{Amb05},
the use of the  Alberti rank one theorem can be circumvented. In our
proof, it is possible to use the following, much simpler ingredient. For any
matrix $M$, we have
\[
\inf_\rho \int |\langle M z,\nabla \rho(z) \rangle | \,dz = |\tr(M)|\, ,
\]
where the infimum is taken over all smooth kernels with total mass
one. Essentially applying this result to the matrices $D^s b(x)$ we
manipulate in the proof, and using  $\tr(D^s b(x)) = \diver(b^s)(x)=0$,
we may obtain an estimate analogous to \eqref{lastbound}, thus
\eqref{lastbound0}. We then conclude our
argument similarly. We however believe that considering the Alberti rank one Theorem
helps to better understand the geometry of the problem.
\end{remark}



{\bf Acknowledgements.}\\
This work was mostly completed while the second author was visiting the
Institute for Mathematics and its Applications and the Department of Mathematics
of the University of Minnesota. The hospitality of these institutions  is
gratefully acknowledged. Both authors would also like to thank Pierre-Emmanuel
Jabin for helpful discussions.

\bibliographystyle{plain}
{\bibliography{BVflow_HLB}}

\end{document}